\documentclass[11pt]{article}

\usepackage{amsmath,amssymb}
\usepackage{latexsym}
\usepackage{graphicx}
\usepackage{bm}

\newtheorem{thm}{Theorem}%[section]

\newenvironment{proof}{\begin{trivlist}
                       \item[]{\bf Proof.}
                       \hspace{0cm}}{\hfill $\Box$
                       \end{trivlist}}

\begin{document}
\title{Dynamical systems method for solving linear\\
finite-rank operator equations}

\author{N. S. Hoang$\dag$\footnotemark[1]
\quad 	  A. G. Ramm$\dag$\footnotemark[3] \\
\\
\\
$\dag$Mathematics Department, Kansas State University,\\
Manhattan, KS 66506-2602, USA
}

\renewcommand{\thefootnote}{\fnsymbol{footnote}}
\footnotetext[1]{Email: nguyenhs@math.ksu.edu}
\footnotetext[3]{Corresponding author. Email: ramm@math.ksu.edu}

\date{}
\maketitle

\begin{abstract}
\noindent
A version of the Dynamical Systems Method (DSM) for solving ill-conditioned linear algebraic systems 
is studied in this paper. 
An {\it a priori} and {\it a posteriori} stopping rules are justified. 
An iterative scheme is constructed for solving ill-conditioned linear algebraic systems.

\end{abstract}

\section{Introduction}

We want to solve stably the equation \begin{equation} \label{eq1} Au=f,
\end{equation} where A is a linear bounded operator in a real Hilbert
space $H$. We assume that \eqref{eq1} has a solution, possibly nonunique,
and denote by $y$ the unique minimal-norm solution to \eqref{eq1}, $y\perp
\mathcal{N}:=\mathcal{N}(A):= \{u: Au=0\}$, $Ay=f$. We assume that the
range of $A$, $R(A)$, is not closed, so problem \eqref{eq1} is ill-posed.
Let $f_\delta$, $\|f-f_\delta\|\leq\delta$, be the noisy data. We want to
construct a stable approximation of $y$, given $\{\delta, f_\delta, A\}$.
There are many methods for doing this, see, e.g., 
\cite{I}--\cite{M}, \cite{R499}, \cite{VV}, \cite{V}, to mention some
(of the many) books, where
variational regularization, quasisolutions, quasiinversion, and iterative
regularization are studied, and \cite{R499}-\cite{R491}, where the 
Dynamical Systems 
Method (DSM) is studied systematically (see also \cite{R401}, \cite{VV}, 
\cite{T},
and references therein for related results). 
The
basic new results of this paper are: 1) a new version of the DSM for
solving equation \eqref{eq1} is justified; 2) a stable method for solving
equation \eqref{eq1} with noisy data by the DSM is given; a priori and a
posteriori stopping rules are proposed and justified; 3) an iterative
method for solving linear ill-conditioned algebraic systems, based on the
proposed version of DSM, is formulated; its convergence is proved; 4)
numerical results are given; these results show that the proposed method
yields a good alternative to some of the standard methods (e.g., to
variational regularization, Landweber iterations, and some other methods).

The DSM version we study in this paper consists of solving the Cauchy problem
\begin{equation}
\label{eqi2}
\dot{u}(t)=- P(Au(t)-f),\quad u(0)=u_0, \quad u_0\perp \mathcal{N},\quad \dot{u}:=\frac{du}{dt},
\end{equation}
and proving the existence of the limit
$\lim_{t\to\infty}u(t)=u(\infty)$, and the relation $u(\infty)=y$, i.e.,
\begin{equation}
\label{eqi3}
\lim_{t\to\infty}\|u(t)-y\|=0.
\end{equation}
Here $P$ is a bounded operator such that $T:=PA\ge 0$ is selfadjoint, $\mathcal{N}(T)=\mathcal{N}(A)$.

For any linear (not necessarily bounded) operator $A$ there exists a bounded operator $P$ such
that $T=PA\ge 0$. For example, if $A=U|A|$ is the polar decomposition of $A$,
then $|A|:=(A^*A)^{\frac{1}{2}}$ is a selfadjoint operator, $T:=|A|\ge 0$,
$U$ is a partial isometry, $\|U\|=1$, and if $P:=U^*$, then $\|P\|=1$ and $PA=T$.
Another choice of $P$, namely, $P=(A^*A + aI)^{-1}A^*$, $a=const>0$, is used in Section~\ref{numsec}.  

If the noisy data $f_\delta$ are given, $\|f_\delta-f\|\le \delta$, then we solve the problem
\begin{equation}
\label{eqi4}
\dot{u}_\delta(t)=-P(Au_\delta(t)-f_\delta),\quad u_\delta(0)=u_0, 
\end{equation}
and prove that, for a suitable stopping time $t_\delta$, and $u_\delta :=u_\delta(t_\delta)$,
one has
\begin{equation}
\label{eqi5}
\lim_{\delta\to 0}\|u_\delta-y\|=0.
\end{equation}
An {\it a priori} and an {\it a posteriori} methods for choosing $t_\delta$ are given.

In Section 2 these results are formulated and recipes for choosing $t_\delta$ are proposed.
In Section \ref{numsec} a numerical example is presented.

\section{Formulation and results}

Suppose $ A: H\to H$ is a linear bounded operator in a real Hilbert space $H$.
Assume that equation \eqref{eq1}
has a solution not necessarily unique. Denote by $y$ 
the unique minimal-norm solution i.e., $y\perp \mathcal{N}:=\mathcal{N}(A)$. 
Consider the DSM \eqref{eqi2}
%\begin{equation}
%\label{eq2}
%\begin{split}
%\dot{u}&=-P(Au - f),\\
%u(0)&=u_0,
%\end{split}
%\end{equation}
where $u_0\perp \mathcal{N}$ is arbitrary. Denote 
\begin{equation}
\label{eq6x}
T:=PA,\quad Q:=AP.
\end{equation}
The unique solution to \eqref{eqi2} is
\begin{equation}
\label{2eq7}
u(t)=e^{-tT}u_0 + e^{-tT}\int_0^t e^{sT} ds Pf.
\end{equation}
Let us first show that any ill-posed linear equation \eqref{eq1} with exact data can be solved by the DSM.

\subsection{Exact data}

The following result is known (see \cite{R499}) but a short proof is included for completeness.
\begin{thm}
\label{thm1}
Suppose $u_0\perp \mathcal{N}$ and $T^*=T\ge 0$. Then problem \eqref{eqi2} has a unique solution defined on
$[0,\infty)$, and $u(\infty)=y$, where $u(\infty)=\lim_{t\to\infty} u(t)$.
\end{thm}

\begin{proof}
Denote $w:=u(t)-y,\, w_0:=w(0)=u_0-y$. Note that $w_0\perp \mathcal{N}$. One has
\begin{equation}
\label{eq3}
\dot{w}=-Tw, \quad T:=PA,\quad w(0)=u_0 - y.
\end{equation}
The unique solution to \eqref{eq3} is
$w=e^{-tT}w_0$. Thus,
$$
\|w\|^2=\int_0^{\|T\|} e^{-2t\lambda}d\langle E_\lambda w_0,w_0\rangle.
$$ 
where $\langle u,v\rangle$ is the inner product in $H$, and $E_\lambda$ is 
the resolution of the
identity of $T$. Thus,
$$
\|w(\infty)\|^2 =\lim_{t\to \infty}\int_0^{\|T\|} e^{-2t\lambda}d\langle E_\lambda w_0,w_0\rangle=\|P_\mathcal{N}w_0\|^2=0,
$$
where $P_\mathcal{N}=E_0-E_{-0}$ is the orthogonal projector onto $\mathcal{N}$.
Theorem \ref{thm1} is proved.
\end{proof}

\subsection{Noisy data $f_\delta$}
\label{Sec2.2}

Let us solve stably equation \eqref{eq1} assuming that $f$ is not known,
but $f_\delta$, the noisy data, are known, 
where $\|f_\delta-f\|\le \delta$. Consider the following 
DSM 
\begin{equation}
\label{eq7delta}
\dot{u}_\delta = - P(Au_\delta - f_\delta),\quad u_\delta(0)=u_0.
\end{equation}
Denote 
$$
w_\delta:=u_\delta-y,\quad T:=PA,\quad w_\delta(0)=w_0:=u_0-y\in \mathcal{N}^\perp.
$$ 
Let us prove the following result:

\begin{thm}
\label{thm2}
If $T=T^*\ge0$, $\lim_{\delta\to0}t_\delta = \infty,\,\lim_{\delta\to 0}t_\delta \delta=0$, and
$w_0\in \mathcal{N}^\perp$, 
then 
$$
\lim_{\delta\to 0}\|w_\delta(t_\delta)\|=0.
$$
\end{thm}

\begin{proof}
One has
\begin{equation}
\label{eq4}
\dot{w}_\delta= -Tw_\delta + \zeta_\delta,\quad \, \zeta_\delta
=P(f_\delta - f),\quad \|\zeta_\delta\|\le\|P\|\delta.
\end{equation}
The unique solution of equation \eqref{eq4} is
$$
w_\delta(t)=e^{-tT}w_\delta(0)+\int_0^te^{-(t-s)T}\zeta_\delta ds.
$$
Let us show that 
$\lim_{\delta\to 0} \|w_\delta(t_\delta)\|=0$. 
One has
\begin{equation}
\label{extra1}
\lim_{t\to\infty} \|w_\delta(t)\| \le \lim_{t\to\infty}\|e^{-tT}w_\delta(0)\|
+\lim_{t\to\infty}\bigg{\|}\int_0^te^{-(t-s)T}\zeta_\delta ds\bigg{\|}.
\end{equation}
Let $E_\lambda$ be the resolution of identity corresponding to $T$. 
One uses the spectral theorem and gets:
\begin{equation}
\label{eq5}
\begin{split}
\int_0^te^{-(t-s)T}ds\zeta_\delta&=\int_0^t\int_0^{\|T\|} dE_\lambda \zeta_\delta e^{-(t-s)\lambda} ds\\
&=\int_0^{\|T\|} e^{-t\lambda}\frac{e^{t\lambda}-1}{\lambda}dE_\lambda\zeta_\delta
=\int_0^{\|T\|}\frac{1-e^{-t\lambda}}{\lambda}dE_\lambda\zeta_\delta.
\end{split}
\end{equation}
Note that
\begin{equation}
\label{eq6}
0\le\frac{1-e^{-t\lambda}}{\lambda}\le t,\quad \forall \lambda>0,\quad t\ge 0,
\end{equation}
since $1-x\le e^{-x}$ for $x\ge 0$. 
From \eqref{eq5} and \eqref{eq6}, one obtains
\begin{equation}
\label{extra2}
\begin{split}
\bigg{\|}\int_0^te^{-(t-s)T}ds\zeta_\delta\bigg{\|}^2
&=\int_0^{\|T\|}\big{|}\frac{1-e^{-t\lambda}}{\lambda}\big{|}^2d\langle E_\lambda\zeta_\delta,\zeta_\delta\rangle\\
&\le t^2\int_0^{\|T\|} d\langle E_\lambda\zeta_\delta,\zeta_\delta\rangle\\
&=t^2\|\zeta_\delta\|^2.
\end{split}
\end{equation}
Since $\|\zeta_\delta\|\le \|P\|\delta$, from \eqref{extra1} and \eqref{extra2}, one gets
$$
\lim_{\delta\to0} \|w_\delta(t_\delta)\| \le \lim_{\delta\to 0}\bigg{(}
\| e^{-t_\delta T}w_\delta(0)\|+t_{\delta}\delta\|P\|\bigg{)}=0.
$$
Here we have used the relation:
$$
\lim_{\delta\to 0}\| e^{-t_\delta T}w_\delta(0)\|=\|P_\mathcal{N}w_0\|=0,
$$
and the last equality holds because $w_0\in \mathcal{N}^\perp$. 
Theorem \ref{thm2} is proved.
\end{proof}

From Theorem \ref{thm2}, it follows that the relation 
$$
t_\delta=\frac{C}{\delta^\gamma},\quad \gamma=\text{const},\quad \gamma\in(0,1)
$$ 
where $C>0$ is a constant, 
can be used as an \textit{a priori} stopping rule, i.e., for such $t_\delta$ one has
\begin{equation}\
\label{eq7}
\lim_{\delta\to0}\|u_\delta(t_\delta)-y\|=0.
\end{equation}

\subsection{Discrepancy principle}

In this section we assume that $A$ is a linear finite-rank operator. Thus, it is
a linear bounded operator. 
Let us consider equation \eqref{eq1} with noisy data $f_\delta$, and a DSM of the form
\begin{equation}
\label{eq8}
\dot{u}_\delta = - PA u_\delta + Pf_\delta,\quad u_\delta(0)=u_0.
\end{equation}
for solving this equation. 
Equation \eqref{eq8} has been used in Section~\ref{Sec2.2}. 
Recall that $y$ denotes the minimal-norm solution of equation \eqref{eq1}.
Example of a choice of $P$ is given in Section~\ref{numsec}.

\begin{thm}
\label{thm3}
Let $T:=PA$, $Q:=AP$.
Assume that $\|Au_0-f_\delta\|> C\delta$, $Q=Q^*\ge 0$, $T^*=T\ge0$, $T$ is a finite-rank operator. Let
$\mathcal{N}(T)=:\mathcal{N}$. Note that $\mathcal{N}(T)=\mathcal{N}(A)$. The
solution $t_\delta$ to the equation
\begin{equation}
\label{eq9}
h(t):=\|Au_\delta(t)- f_\delta\|=C\delta,\quad C=\text{const},\quad C\in (1,2),
\end{equation}
does exist, is unique, and
\begin{equation}
\label{eq10}
\lim_{\delta\to 0} \|u_\delta(t_\delta)-y\|=0,
\end{equation}
where $y$ is the unique minimal-norm solution to \eqref{eq1}.
\end{thm}

\begin{proof}
Denote 
$$
v_\delta(t):=Au_\delta(t)- f_\delta,\quad w(t):=u(t)-y,\quad w_0:=u_0-y.
$$ 
One has
\begin{equation}
\label{eq11}
\begin{split}
\frac{d}{dt}\|v_\delta(t)\|^2
&= 2\langle A\dot{u}_\delta(t),Au_\delta(t)-f_\delta \rangle\\
&= 2\langle A[-P(Au_\delta(t) - f_\delta)],Au_\delta(t)-f_\delta \rangle\\
&=-2\langle AP(Au_\delta-f_\delta),Au_\delta-f_\delta\rangle\le 0.
\end{split}
\end{equation}
where the last inequality holds because $AP=Q\ge0$.
Thus, $\|v_\delta(t)\|$ is a nonincreasing function.

Let us prove that equation \eqref{eq9} has a solution for $C\in (1,2)$. 
One has the following commutation formulas: 
$$
e^{-sT}P=Pe^{-sQ},\quad Ae^{-sT}=e^{-sQ}A.
$$
Using these formulas and the representation 
$$
u_\delta(t)=e^{-tT}u_0+\int_0^te^{-(t-s)T}Pf_\delta ds,
$$ 
one gets:
\begin{equation}
\label{eq0302}
\begin{split}
v_\delta(t)
&= Au_\delta(t)-f_\delta\\
&= Ae^{-tT}u_0+A\int_0^te^{-(t-s)T}Pf_\delta ds -f_\delta \\
&= e^{-t Q}Au_0+e^{-t Q}\int_0^{t} e^{sQ}dsQf_\delta-f_\delta \\
%&= e^{-tQ}Au_0 - e^{-tQ}f_\delta\\
&= e^{-tQ}A(u_0-y)+e^{-tQ}f+e^{-tQ}(e^{tQ}-I)f_\delta - f_\delta\\
&= e^{-t Q}Aw_0 -e^{-t Q}f_\delta +e^{-t Q}f= e^{-tQ}Au_0 - e^{-tQ}f_\delta.
\end{split}
\end{equation}
%where $Q:=AA^*$. Let $w(t):=u(t)-y,\, w_0=u_0-y$. 
Note that 
$$
\lim_{t\to\infty}e^{-t Q}Aw_0=\lim_{t\to\infty}Ae^{-t T}w_0 = AP_\mathcal{N}w_0=0.
$$ 
Here the continuity of $A$ and the following
relation
$$
\lim_{t\to\infty}e^{-tT}w_0=\lim_{t\to\infty}\int_0^{\|T\|}e^{-st}dE_sw_0=(E_0-E_{-0})w_0=P_\mathcal{N}w_0,
$$
were used. 
Therefore,
\begin{equation}
\label{eq12}
\lim_{t\to\infty}\|v_\delta(t)\|=\lim_{t\to\infty}\|e^{-t Q}(f-f_\delta)\|\le
\|f-f_\delta\|\le\delta,
\end{equation}
where $\|e^{-tQ}\|\le 1$ because $Q\ge0$. The function $h(t)$ is continuous on $[0,\infty)$, 
$h(0)=\|Au_0-f_\delta\|>C\delta$, $h(\infty)\le \delta$.
Thus, equation \eqref{eq9} must have a solution $t_\delta$.

Let us prove the uniqueness of $t_\delta$. If $t_\delta$ is non-unique, then without loss of generality we can assume that there exists 
$t_1>t_\delta$ such that $\|Au_\delta(t_1)- f_\delta\|=C\delta$. Since $\|v_\delta(t)\|$ is 
nonincreasing and $\|v_\delta(t_\delta)\|=\|v_\delta(t_1)\|$, one has
$$
\|v_\delta(t)\|=\|v_\delta(t_\delta)\|,\quad \forall t\in [t_\delta, t_1].
$$
Thus,
\begin{equation}
\label{eq13}
\frac{d}{dt}\|v_\delta(t)\|^2=0,\quad \forall t\in (t_\delta, t_1).
\end{equation}
Using \eqref{eq11} and \eqref{eq13} one obtains
$$
\| \sqrt{AP}(Au_\delta(t)-f_\delta)\|^2=\langle AP(Au_\delta(t)-f_\delta), Au_\delta(t)-f_\delta\rangle= 0,\quad \forall t\in [t_\delta,t_1],
$$
where $\sqrt{AP}=Q^{\frac{1}{2}}\ge 0$ is well defined since $Q=Q^*\ge 0$.
This implies $Q^{\frac{1}{2}}(Au_\delta-f_\delta)=0$. Thus
\begin{equation}
\label{thich1}
\begin{split}
Q(Au_\delta(t)-f_\delta)=0,\quad \forall t\in [t_\delta,t_1].
\end{split}
\end{equation}
From \eqref{eq0302} one gets:
\begin{equation}
\label{3eq23}
v_\delta(t)= Au_\delta(t)-f_\delta= e^{-tQ}Au_0 - e^{-tQ}f_\delta.
\end{equation}
Since $Qe^{-tQ}=e^{-tQ}Q$ and $e^{-tQ}$ is an isomorphism, equalities \eqref{thich1} and \eqref{3eq23} imply
\begin{align*}
 Q(Au_0 - f_\delta)=0.
\end{align*}
This and \eqref{3eq23} imply 
$$
AP(Au_\delta(t)-f_\delta)=e^{-tQ}(QAu_0 - Qf_\delta)=0,\quad t\ge0.
$$
This and \eqref{eq11} imply
\begin{equation}
\label{eq14}
\frac{d}{dt}\|v_\delta\|^2=0,\quad t\ge0.
\end{equation}
Consequently, 
$$
C\delta<\|Au_\delta(0)-f_\delta\|=\|v_\delta(0)\|
=\|v_\delta(t_\delta)\|
=\|Au_\delta(t_\delta)-f_\delta\|
=C\delta.
$$ 
This is a contradiction which proves
the uniqueness of $t_\delta$.

Let us prove \eqref{eq10}. First, we have the following estimate:
\begin{equation}
\label{eq16}
\begin{split}
\|Au(t_\delta)-f\|&\le \|Au(t_\delta)-Au_\delta(t_\delta)
\|+\|Au_\delta(t_\delta)-f_\delta\|+\|f_\delta -f\|\\
&\le \bigg{\|}e^{-t_\delta Q}\int_0^{t_\delta}e^{sQ}Qds \bigg{\|} \|f_\delta-f\|+C\delta+\delta,
\end{split}
\end{equation}
where $u(t)$ solves \eqref{eqi2} and $u_\delta(t)$ solves \eqref{eq7delta}. 
One uses the inequality:
$$
\big{\|}e^{-t_\delta Q}\int_0^{t_\delta}e^{sQ}Qds \big{\|}=
\|I-e^{-t_\delta Q}\|\le 2,
$$
and concludes from \eqref{eq16}, that
\begin{equation}
\label{eq17}
\lim_{\delta\to0}\|Au(t_\delta)-f\|=0.
\end{equation}
Secondly, 
we claim that 
$$
\lim_{\delta\to0}t_\delta=\infty.
$$ 
Assume the contrary. Then there exist $t_0>0$ and a sequence
$(t_{\delta_n})_{n=1}^\infty$,
$t_{\delta_n}<t_0$, such that
\begin{equation}
\label{eq18}
\lim_{n\to\infty}\|Au(t_{\delta_n})-f\|=0.
\end{equation}
Analogously to \eqref{eq11}, one proves that 
$$
\frac{d}{dt}\|v\|^2\le 0,
$$ 
where $v(t):=Au(t)-f$.
Thus, $\|v(t)\|$ is nonincreasing. 
This and \eqref{eq18} imply the relation $\|v(t_0)\|=\|Au(t_0)-f\|=0$.
Thus,
$$
0=v(t_0)=e^{-t_0Q}A(u_0-y).
$$
This implies $A(u_0-y)=e^{t_0Q}e^{-t_0Q}A(u_0-y)=0$, so $u_0-y\in \mathcal{N}$. 
Since $u_0-y\in \mathcal{N}^\perp$, it follows that
 $u_0=y$. This is a contradiction because 
$$
C\delta\le\|Au_0-f_\delta\|=\|f-f_\delta\|\le\delta, \quad 1<C<2.
$$ 
Thus, 
\begin{equation}
\label{eq03new1}
\lim_{\delta\to0}t_\delta=\infty.
\end{equation}

Let us continue the proof of \eqref{eq10}. 
%If $T$ is a finite rank operator, then we can argue as follows:
From \eqref{eq0302} and the relation $\|Au_\delta(t_\delta)-f_\delta\|=C\delta$, one has
\begin{equation}
\label{5eq30}
\begin{split}
C\delta t_\delta
&=\|t_\delta e^{-t_\delta Q}Aw_0 - t_\delta e^{-t_\delta Q}(f_\delta - f)\|\\
&\le\|t_\delta e^{-t_\delta Q}Aw_0\| + \| t_\delta e^{-t_\delta Q}(f_\delta - f)\|\\
&\le \|t_\delta e^{-t_\delta Q}Aw_0\| + t_\delta \delta. 
\end{split}
\end{equation}
We claim that
\begin{equation}
\label{hoithay}
\lim_{\delta \to 0} t_\delta e^{-t_\delta Q}Aw_0 =\lim_{\delta\to 0}t_\delta Ae^{-t_\delta T}w_0= 0.
\end{equation}
Note that \eqref{hoithay} holds if $T\ge0$ has finite rank, and $w_0\in \mathcal{N}^\perp$. 
It also holds if $T\ge 0$ is compact and the Fourier coefficients 
$w_{0j}:=\langle w_0,\phi_j\rangle$, $T\phi_j=\lambda_j\phi_j$, decay sufficiently fast. 
In this case
$$
\|Ae^{-tT}w_0\|^2\le\|T^{\frac{1}{2}}e^{-tT}w_0\|^2
=\sum_{j=1}^\infty \lambda_je^{-2\lambda_jt}|w_{0j}|^2:=S=o(\frac{1}{t^2}), 
\quad t\to\infty,
$$ 
provided that $\sum_{j=1}^\infty|w_{0j}|\lambda_j^{-2}<\infty$. 
Indeed $S=\sum_{\lambda_j\le\frac{1}{t^{\frac{2}{3}}}}
+\sum_{\lambda_j>\frac{1}{t^{\frac{2}{3}}}}:=S_1+S_2$. 
One has
$$
S_1\le\frac{1}{t^2}\sum_{\lambda_j\le t^{-\frac{2}{3}}}\frac{|w_{0j}|^2}{\lambda_j^2}=o(\frac{1}{t^2}),
\quad S_2\le ce^{-2t^{\frac{1}{3}}}=o(\frac{1}{t^2}),\quad t\to\infty, 
$$
where $c>0$ is a constant.

From \eqref{hoithay} and \eqref{5eq30}, one gets
$$
0\le \lim_{\delta\to0} (C-1)\delta t_\delta \le \lim_{\delta\to 0} \|t_\delta e^{-t_\delta Q}Aw_0\|=0.
$$
Thus, 
\begin{equation}
\label{eq03new2}
\lim_{\delta\to 0}\delta t_\delta=0
\end{equation}
Now \eqref{eq10} follows from \eqref{eq03new1}, 
\eqref{eq03new2} and Theorem~\ref{thm2}. Theorem 3 is proved.
\end{proof}

\subsection{An iterative scheme}
\label{itersec}

Let us solve stably equation \eqref{eq1} assuming that $f$ is not known,
but $f_\delta$, the noisy data, are known, 
where $\|f_\delta-f\|\le \delta$. Consider the following
dicrete version of the DSM: 
\begin{equation}
\label{eq7ndelta}
u_{n+1,\delta} = u_{n,\delta} - hP(Au_{n,\delta} - f_\delta),\quad u_{\delta,0}=u_0.
\end{equation}
Let us denote $u_n:=u_{n,\delta}$ when $\delta\not=0$, and 
set 
$$
w_n:=u_n-y,\quad T:=PA,\quad w_0:=u_0-y\in \mathcal{N}^\perp.
$$ 
Let $n=n_\delta$  be the stopping rule for iterations \eqref{eq7ndelta}.
Let us prove the following result:

\begin{thm}
\label{4thm2}
Assume that $T=T^*\ge0$, $h\|T\|< 2$, 
$\lim_{\delta\to0}n_\delta h = \infty,\,\lim_{\delta\to 0}n_\delta h \delta=0$, 
and $w_0\in \mathcal{N}^\perp$. Then 
\begin{equation}
\label{eq32x}
\lim_{\delta\to 0}\|w_{n_\delta}\|= \lim_{\delta\to 0}\|u_{n_\delta}-y\| =0.
\end{equation}
\end{thm}

\begin{proof}
One has
\begin{equation}
\label{4eq4}
w_{n+1} = w_n -h Tw_n + h \zeta_\delta,\quad \, \zeta_\delta
=P(f_\delta - f),\quad \|\zeta_\delta\|\le\|P\|\delta,\quad w_0=u_0-y.
\end{equation}
The unique solution of equation \eqref{4eq4} is
$$
w_{n+1} = (I-hT)^{n+1}w_0 + h\sum_{i=0}^n(I-hT)^i \zeta_\delta.
$$
Let us show that 
$\lim_{\delta\to 0} \|w_{n_\delta}\|=0$. 
One has
\begin{equation}
\label{4extra1}
\|w_n\| \le \|(I-hT)^{n}w_0\|
+ \bigg{\|}h\sum_{i=0}^{n-1}(I-hT)^i \zeta_\delta\bigg{\|}.
\end{equation}
Let $E_\lambda$ be the resolution of identity corresponding to $T$. 
One uses the spectral theorem and gets:
\begin{equation}
\label{4eq5}
\begin{split}
h\sum_{i=0}^{n-1}(I-hT)^i &= h\sum_{i=0}^{n-1}\int_0^{\|T\|} (1-h\lambda)^i dE_\lambda\\
&=h\int_0^{\|T\|} \frac{1 - (1-\lambda h)^{n}}{1-(1-h\lambda)}dE_\lambda
=\int_0^{\|T\|} \frac{1 - (1-\lambda h)^{n}}{\lambda}dE_\lambda.
\end{split}
\end{equation}
Note that
\begin{equation}
\label{4eq6}
0\le\frac{1-(1-h\lambda )^{n}}{\lambda}\le hn,\quad \forall \lambda>0,\quad t\ge 0,
\end{equation}
since $1-(1-\alpha)^n\le \alpha n$ for all $\alpha \in [0,2]$. 
From \eqref{4eq5} and \eqref{4eq6}, one obtains
\begin{equation}
\label{4extra2}
\begin{split}
\bigg{\|} h\sum_{i=0}^{n-1}(I-hT)^i\zeta_\delta\bigg{\|}^2
&=\int_0^{\|T\|}\big{|}\frac{1 - (1-\lambda h)^{n}}{\lambda}\big{|}^2d\langle E_\lambda\zeta_\delta,\zeta_\delta\rangle\\
&\le (hn)^2\int_0^{\|T\|} d\langle E_\lambda\zeta_\delta,\zeta_\delta\rangle\\
&= (nh)^2\|\zeta_\delta\|^2.
\end{split}
\end{equation}
Since $\|\zeta_\delta\|\le \|P\|\delta$, from \eqref{4extra1} and \eqref{4extra2}, one gets
$$
\lim_{\delta\to0} \|w_{n_\delta}\| \le \lim_{\delta\to 0}\bigg{(}
\| (I-hT)^{n_\delta}w_\delta(0)\|+ hn_\delta \delta\|P\|\bigg{)}=0.
$$
Here we have used the relation:
$$
\lim_{\delta\to 0}\|(I-hT)^{n_\delta}w_\delta(0)\|=\|P_\mathcal{N}w_0\|=0,
$$
and the last equality holds because $w_0\in \mathcal{N}^\perp$. 
Theorem \ref{4thm2} is proved.
\end{proof}

From Theorem \ref{4thm2}, it follows that the relation 
$$
n_\delta=\frac{C}{h\delta^\gamma},\quad \gamma=\text{const},\quad \gamma\in(0,1)
$$ 
where $C>0$ is a constant, 
can be used as an \textit{a priori} stopping rule, i.e., for such $n_\delta$ one has
\begin{equation}\
\label{4eq7}
\lim_{\delta\to0}\|u_{n_\delta}-y\|=0.
\end{equation}

\subsection{An iterative scheme with a stopping rule 
based on a discrepancy principle}

In this section we assume that $A$ is a linear finite-rank operator. 
Thus, it is
a linear bounded operator. 
Let us consider equation \eqref{eq1} with noisy data $f_\delta$, and a DSM of the form
\begin{equation}
\label{5eq8}
u_{n+1} = u_n - hP(A u_n - f_\delta),\quad u_0 =u_0,
\end{equation}
for solving this equation. 
Equation \eqref{5eq8} has been used in Section~\ref{itersec}. 
Recall that $y$ denotes the minimal-norm solution of equation \eqref{eq1}.
Example of a choice of $P$ is given in Section~\ref{numsec}.

Note that $\mathcal{N}:=\mathcal{N}(T)=\mathcal{N}(A)$.

\begin{thm}
\label{5thm3}
Let $T:=PA$, $Q:=AP$.
Assume that $\|Au_0-f_\delta\|> C\delta$, $Q=Q^*\ge 0$, $T^*=T\ge0$, 
$h\|T\|< 2$,  $h\|Q\|< 2$, and
$T$ is a finite-rank operator. 
Then there exists a unique $n_\delta$
such that
\begin{equation}
\label{5eq9}
\|Au_{n_\delta}- f_\delta\|\le C\delta < \|Au_{n_\delta-1}- f_\delta\|,\quad C=\text{const},\quad C\in (1,2).
\end{equation}
For this $n_\delta$ one has:
\begin{equation}
\label{5eq10}
\lim_{\delta\to 0} \|u_{n_\delta} -y\|=0.
\end{equation}
\end{thm}

\begin{proof}
Denote 
$$
v_n:=Au_n - f_\delta,\quad w_n:=u_n-y,\quad w_0:=u_0-y.
$$ 
From \eqref{5eq8}, one gets
\begin{align*}
v_{n+1} &= Au_{n+1} - f_\delta = Au_n -f_\delta - h AP(Au_n - f_\delta)  
= v_n - h Qv_n.  
\end{align*}
This implies
\begin{equation}
\label{5eq11}
\begin{split}
\|v_{n+1}\|^2 - \|v_{n}\|^2
&= \langle v_{n+1} - v_n, v_{n+1} + v_n \rangle\\ 
&= \langle -hQ v_n, v_n - hQv_n + v_n \rangle\\
&= -\langle v_n, hQ(2-hQ)v_n  \rangle\le 0\\
\end{split}
\end{equation}
where the last inequality holds because $AP=Q\ge0$ and $\|hQ\| < 2$.
Thus, $(\|v_n\|)_{n=1}^\infty$ is a nonincreasing sequence.

Let us prove that equation \eqref{5eq9} has a solution for $C\in (1,2)$. 
One has the following commutation formulas: 
$$
(I - hT)^nP=P(I - hQ)^n,\quad A(I - hT)^n=(I - hQ)^n A.
$$
Using these formulas, the representation 
$$
u_n = (I - hT)^nu_0 + h\sum_{i=0}^{n-1}(I - hT)^iPf_\delta,
$$ 
and the identity $(I-B)\sum_{i=0}^{n-1}B^i=I-B^n$, with $B=I-hQ$, $I-B=hQ$, 
one gets:
\begin{equation}
\label{5eq0302}
\begin{split}
v_n
&= Au_n-f_\delta\\
&= A(I - hT)^n u_0 +Ah\sum_{i=0}^{n-1}(I - hT)^i Pf_\delta -f_\delta\\
&= (I - hQ)^nAu_0+\sum_{i=0}^{n-1}(I - hQ)^ihQf_\delta-f_\delta \\
&= (I - hQ)^nAu_0 - (I-(I - hQ)^n)f_\delta-f_\delta\\
&= (I - hQ)^n(Au_0-f)+(I - hQ)^n(f- f_\delta)\\
&= (I - hQ)^nAw_0 + (I - hQ)^n(f-f_\delta).%\\
\end{split}
\end{equation}
If $V=V^*\geq 0$ is an operator with $||V||\leq 2$, then
$||I-V||=\sup_{0\leq s \leq 2}|1-s|\leq 1$.
 
Note that 
$$
\lim_{n\to\infty}(I - hQ)^n Aw_0=\lim_{n\to\infty}A(I - hT)^nw_0 = 
AP_\mathcal{N}w_0=0,
$$ 
where $P_\mathcal{N}$ is the orthoprojection onto the null-space
$\mathcal{N}$ of the operator $T$, and
the continuity of $A$ and the following
relation
$$
\lim_{n\to\infty}(I - hT)^n w_0=\lim_{n\to\infty}\int_0^{\|T\|}(1-s 
h)^ndE_sw_0
=(E_0-E_{-0})w_0=P_\mathcal{N}w_0,\quad 0\le sh < 2,
$$
were used. 
Therefore,
\begin{equation}
\label{5eq12}
\lim_{n\to\infty}\|v_\delta(t)\|=\lim_{n\to\infty}\|(I-hQ)^n(f-f_\delta)\|\le
\|f-f_\delta\|\le\delta,
\end{equation}
where $\|I-hQ\|\le 1$ because $Q\ge0$ and $\|hQ\|<2$. The sequence 
$\{\|v_n\|\}_{n=1}^\infty$ is nonincreasing with 
$\|v_0\|>C\delta$ and $\lim_{n\to\infty}\|v_n\|\le\delta$. Thus, there 
exists $n_\delta >0$ such that \eqref{5eq9} holds.

Let us prove \eqref{5eq10}. Let $u_{n,0}$ be the sequence defined by the relations:
$$
u_{n+1,0} = u_{n,0} - hP(Au_{n,0}-f),\quad u_{0,0}=u_0.
$$
First, we have the following estimate:
\begin{equation}
\label{5eq16}
\begin{split}
\|Au_{n_\delta,0}-f\|&\le \|Au_{n_\delta}-Au_{n_\delta,0}\|+\|Au_{n_\delta}-f_\delta\|+\|f_\delta -f\|\\
&\le \bigg{\|}\sum_{i= 0}^{n_\delta -1}(I - hQ)^ihQ\bigg{\|} \|f_\delta-f\|+C\delta+\delta.
\end{split}
\end{equation}
Since $0\le hQ<2$, one has $ ||I-hQ||\leq 1$. This implies the following 
inequality:
$$
\bigg{\|}\sum_{i= 0}^{n_\delta -1}(I - hQ)^ihQ\bigg{\|}=\|I-(I - hQ)^{n_\delta}\|\le 2,
$$
and concludes from \eqref{5eq16}, that
\begin{equation}
\label{5eq17}
\lim_{\delta\to0}\|Au_{n_\delta,0}-f\|=0.
\end{equation}
Secondly, 
we claim that 
$$
\lim_{\delta\to0}hn_\delta=\infty.
$$ 
Assume the contrary. Then there exist $n_0>0$ and a sequence
$(n_{\delta_n})_{n=1}^\infty$,
$n_{\delta_n}<n_0$, such that
\begin{equation}
\label{5eq18}
\lim_{n\to\infty}\|Au_{n_\delta,0}-f\|=0.
\end{equation}
Analogously to \eqref{5eq11}, one proves that 
$$
\|v_{n,0}\|\le\|v_{n-1,0}\|,
$$ 
where $v_{n,0}=Au_{n,0}-f$.
Thus, the sequence $\|v_{n,0}\|$ is nonincreasing. 
This and \eqref{5eq18} imply the relation $\|v_{n_0,0}\|=\|Au_{n_0,0}-f\|=0$.
Thus,
$$
0=v_{n_0,0}=(I - hQ)^{n_0}A(u_0-y).
$$
This implies $A(u_0-y)=(I - hQ)^{-n_0}(I - hQ)^{n_0}A(u_0-y)=0$, so $u_0-y\in \mathcal{N}$. 
Since, by the assumption, $u_0-y\in \mathcal{N}^\perp$, it follows that
 $u_0=y$. This is a contradiction because 
$$
C\delta\le\|Au_0-f_\delta\|=\|f-f_\delta\|\le\delta, \quad 1<C<2.
$$ 
Thus, 
\begin{equation}
\label{5eq03new1}
\lim_{\delta\to0}hn_\delta=\infty.
\end{equation}

Let us continue the proof of \eqref{5eq10}. 
From \eqref{5eq0302} and $\|Au_{n_\delta}-f_\delta\|=C\delta$, one has
\begin{equation}
\label{eq47x}
\begin{split}
C\delta n_\delta h
&=\|n_\delta h(I - hQ)^{n_\delta}Aw_0 - n_\delta h(I - hQ)^{n_\delta}(f_\delta - f)\|\\
&\le\|n_\delta h(I - hQ)^{n_\delta}Aw_0\| + \|n_\delta h(I - hQ)^{n_\delta}(f_\delta - f)\|\\
&\le \|n_\delta h(I - hQ)^{n_\delta}Aw_0\| + n_\delta h\delta. 
\end{split}
\end{equation}
We claim that if $w_0\in\mathcal{N}^\perp$, $0\le hT<2$, and $T$ is a finite-rank operator, then
\begin{equation}
\label{5hoithay}
\lim_{\delta \to 0} n_\delta h(I - hQ)^{n_\delta}Aw_0 =\lim_{\delta\to 0}n_\delta h A(I - hT)^{n_\delta}w_0= 0.
\end{equation}
From \eqref{eq47x} and \eqref{5hoithay} one gets
$$
0\le \lim_{\delta\to0} (C-1)\delta hn_\delta \le \lim_{\delta\to 0} \|n_\delta h(I - hQ)^{n_\delta}Aw_0\|=0.
$$
Thus, 
\begin{equation}
\label{5eq03new2}
\lim_{\delta\to 0}\delta n_\delta h=0
\end{equation}
Now \eqref{5eq10} follows from \eqref{5eq03new1}, \eqref{5eq03new2} and Theorem~\ref{4thm2}.
Theorem~\ref{5thm3} is proved.
\end{proof}

\section{Numerical experiments}
\label{numsec}

\subsection{Computing $u_\delta(t_\delta)$}

In \cite{R540} an DSM \eqref{eq7delta} was investigated with $P=A^*$ and
the SVD of $A$ was assumed known. In general, it is computationally
expensive to get the SVD of large scale matrices. In this paper, we have
derived an iterative scheme for solving ill-conditioned linear algebraic
systems $Au=f_\delta$ without using SVD of $A$.

Choose $P=(A^*A+a)^{-1}A^*$ where $a$ is a fixed positive constant. 
This choice of $P$ satisfies all the conditions in Theorem~\ref{thm3}.
In particular, $Q=AP=A(A^*A+aI)^{-1}A^*=AA^*(AA^* + aI)^{-1}\ge 0$ is 
a selfadjoint operator, and $T=PA=(A^*A+aI)^{-1}A^*A\ge 0$
is a selfadjoint operator.
Since 
$$
\|T\|=\bigg{\|}\int_0^{\|A^*A\|} \frac{\lambda}{\lambda + a}dE_\lambda\bigg{\|}
=\sup_{0\le \lambda\le \|A^*A\|}\frac{\lambda}{\lambda + a}<1,
$$ 
where $E_\lambda$ is the resolution of the identity of $A^*A$, 
the condition
$h\|T\|< 2$ in Theorem~\ref{5thm3} is satisfied for all $ 0<h\le 1$.
Set $h=1$ and $P=(A^*A+a)^{-1}A^*$ in \eqref{5eq8}. Then 
one gets the following iterative scheme:
\begin{equation}
\label{eq30}
u_{n+1} = u_n - (A^*A+aI)^{-1}(A^*Au_n - A^*f_\delta),\quad u_0=0.
\end{equation}
For simplicity we have chosen $u_0=0$. However, one may choose $u_0=v_0$ if
$v_0$ is known to be a better approximation to $y$ than $0$ and $v_0\in 
\mathcal{N}^\perp$.
In iterations \eqref{eq30} we use a stopping rule of discrepancy type. Indeed, we will stop iterations if 
$u_n$ satisfies the following condition
\begin{equation}
\|Au_n - f_\delta\| \le 1.01 \delta.
\end{equation}

The choice of $a$ affects both the accuracy and the computation time of
the method. If $a$ is too large, one needs more iterations to
approach the desired accuracy, so the 
computation time will be large.  
If $a$ is too small then the results become less accurate because
for too small $a$ the inversion of the operator $A^*A+aI$ is an ill-posed 
problem since the operator $A^*A$ is not boundedly invertible.  Using 
the idea
of the choice of the initial guess of regularization parameter in
\cite{R526}, we choose $a$ to satisfy the following condition:

\begin{equation}
\label{eq32}
\delta\le \phi(a):=\|A(A^*A+a)^{-1}A^*f_\delta - f_\delta\| \le 2\delta.
\end{equation}
This can be done by using the following strategy:
\begin{enumerate}
\item{Choose $a:=\frac{\delta \|A\|^2}{3\|f_\delta\|}$ as an initial guess for $a$.}
\item{Compute $\phi(a)$. If $a$ satisfying \eqref{eq32} we are done. Otherwise, we go to step 3.}

\item{If $c=\frac{\phi(a)}{\delta} > 3$ 
we replace $a$ by $\frac{a}{2(c-1)}$ and go back to step 2.
If $2<c\le 3$ then we replace $a$ by $\frac{a}{2(c-1)}$ and go back to step 2.
Otherwise, we go to step 4.}

\item{
If $c=\frac{\phi(a)}{\delta} <1$ %and the inequality $c < 1$ has not occured in iterations 
%then 
we replace $a$ by $3a$. 
If the inequality $c< 1$ has occured in some iteration before, we stop the iteration and use $3a$ 
as our choice for $a$ in iterations \eqref{eq30}}. Otherwise we go back to step 2.
\end{enumerate}

In our experiments, we denote by DSM the iterative scheme \eqref{eq30}, by VR$_i$ 
a Variational Regularization method (VR) with $a$ as the regularization parameter and by VR$_n$ 
the VR in which Newton's method is used for finding the regularization parameter using a 
discrepancy principle. 
We compare these methods in terms of relative error and number of iterations, denoted by 
n$_{iter}$.

All the experiments were carried in double arithmetics
precision environment using MATLAB.

\subsection{A linear algebraic system related to an inverse problem for the heat
equation}

In this section, we apply the DSM and the VR to solve a linear algebraic
system used in \cite{R526}.  This linear algebraic system is a part of
numerical solutions to an inverse problem for the heat equation.  This
problem is reduced to a Volterra integral equation of the first kind with
$[0,1]$ as the integration interval.  The kernel is $K(s,t)=k(s-t)$ with
$$ 
k(t)=\frac{t^{-3/2}}{2\kappa \sqrt{\pi}}\exp(-\frac{1}{4\kappa^2 t}).
$$ 
Here, we use the value $\kappa=1$. In this test in
\cite{R526} the integral equation was discretized by means of simple
collocation and the midpoint rule with $n$ points.  
The unique exact solution
$u_n$ is constructed, and then the right-hand side $b_n$ is produced as
$b_n=A_nu_n$ (see \cite{R526}).  In our test, we use $n=10,20,...,100$
and $b_{n,\delta} = b_n + e_n$, where $e_n$ is a vector containing random
entries, normally distributed with mean 0, variance 1, and scaled so that
$\|e_n\|=\delta_{rel}\|b_n\|$.  This linear system is ill-posed:  
the condition number of $A_{100}$ obtained by using the function
{\it cond} provided in MATLAB is $1.3717\times 10^{37}$.  
%As we have
%discussed earlier, this condition number may not be accurate because
%of the limitations of the program {\it cond} provided in MATLAB.  
%However, 
This number shows that the corresponding linear
algebraic system is severely ill-conditioned.

\begin{table}[ht] 
\caption{Numerical results for the inverse heat equation with $\delta_{rel}=0.05$, $n=10i,\, i=\overline{1,10}$.}
\label{heattab1}
\centering
\small
\begin{tabular}{@{  }c@{\hspace{2mm}}
@{\hspace{2mm}}|c@{\hspace{2mm}}c@{\hspace{2mm}}|c@{\hspace{2mm}}c@{\hspace{2mm}}|
c@{\hspace{2mm}}c@{\hspace{2mm}}|c@{\hspace{2mm}}r@{\hspace{2mm}}l@{}} 
\hline
&\multicolumn{2}{c|}{DSM}&\multicolumn{2}{c|}{VR$_{i}$}&\multicolumn{2}{c|}{VR$_{n}$}\\
$n$&
n$_{\text{iter}}$&$\frac{\|u_\delta-y\|_{2}}{\|y\|_2}$&
n$_{\text{iter}}$&$\frac{\|u_\delta-y\|_{2}}{\|y\|_2}$&
n$_{\text{iter}}$&$\frac{\|u_\delta-y\|_{2}}{\|y\|_2}$\\
\hline
10    &3    &0.1971    &1    &0.2627    &5    &0.2117\\
20    &4    &0.3359    &1    &0.4589    &5    &0.3551\\
30    &4    &0.3729    &1    &0.4969    &5    &0.3843\\
40    &4    &0.3856    &1    &0.5071    &5    &0.3864\\
50    &5    &0.3158    &1    &0.4789    &6    &0.3141\\
60    &6    &0.2892    &1    &0.4909    &6    &0.3060\\
70    &7    &0.2262    &1    &0.4792    &8    &0.2156\\
80    &6    &0.2623    &1    &0.4809    &7    &0.2600\\
90    &5    &0.2856    &1    &0.4816    &7    &0.2715\\
100  & 7    &0.2358    &1    &0.4826    &7    &0.3405\\
%
%10    &2    &0.1296    &1    &0.2081    &5    &0.1389\\
%20    &2    &0.3983    &1    &0.4602    &5    &0.3901\\
%30    &4    &0.2981    &1    &0.4840    &6    &0.2943\\
%40    &4    &0.3283    &1    &0.4828    &5    &0.3432\\
%50    &3    &0.3707    &1    &0.5004    &4    &0.3684\\
%60    &5    &0.2764    &1    &0.4891    &6    &0.2673\\
%70    &5    &0.2805    &1    &0.4822    &7    &0.2779\\
%80    &7    &0.2333    &1    &0.4774    &8    &0.2380\\
%90    &5    &0.2857    &1    &0.4848    &7    &0.2833\\
%100   &5    &0.2986    &1    &0.4917    &6    &0.3017\\
\hline 
\end{tabular}
\end{table}

Table~\ref{heattab1} shows that 
the results obtained by the DSM are comparable to those by the VR$_{n}$ in terms of accuracy. 
The time of computation of the DSM is comparable to that of the VR$_{n}$. 
In some situations, the results by VR$_{n}$ and the DSM are the same although the VR$_n$ uses 3 more iterations than does the DSM.
The conclusion from this Table is that DSM competes favorably with the VR$_{n}$ in both accuracy and time of computation.

\begin{figure}[!h!tb]
\centerline{%
\includegraphics[scale=0.9]{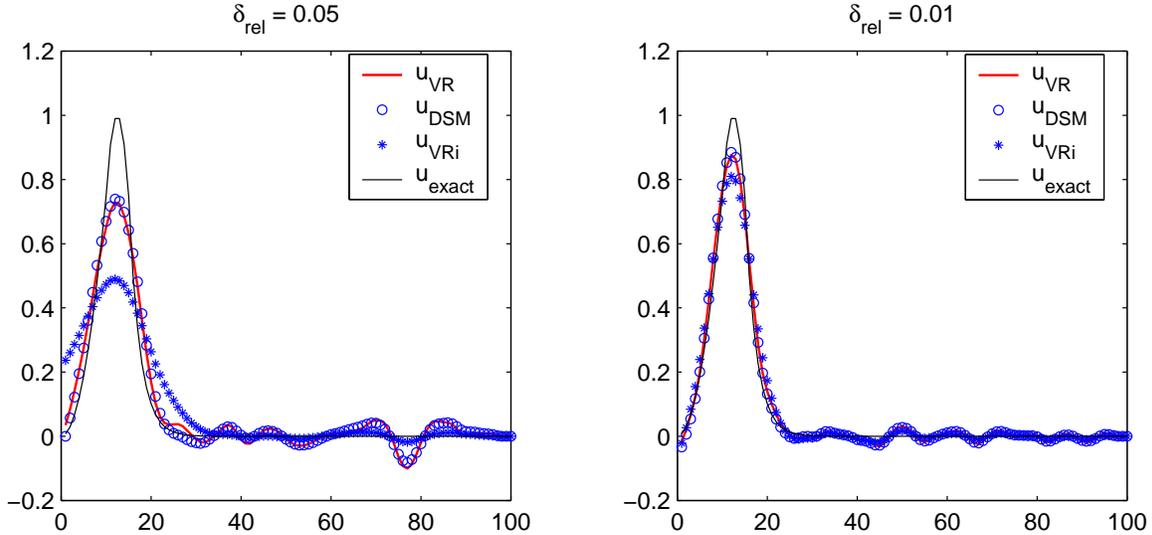}}
%\quad\includegraphics[scale=0.45]{heat2.eps}}
\caption{Plots of solutions obtained by DSM, VR for the inverse heat equation when $n=100$, $\delta_{rel}=0.05$ (left) and $\delta_{rel}=0.01$ (right).}
\label{figheat}
\end{figure}

Figure~\ref{figheat} plots numerical solutions to the inverse heat equation for $\delta_{rel}=0.05$ and $\delta_{rel}=0.01$ when $n=100$. From the figure we can see that the numerical solutions obtained by the DSM are about the same those by the VR$_{n}$. In these examples, the time of computation of the DSM is about the same as that of the VR$_{n}$.

The conclusion is that the DSM competes favorably with the VR$_{n}$ in this experiment.

\section{Concluding remark}

Iterative scheme \eqref{eq30} can be considered as a modification the
Landweber iterations. The difference between the two methods is the multiplication by
$P=(A^*A+aI)^{-1}$.
Our iterative method is much faster than the conventional Landweber iterations.
Iterative method \eqref{eq30} is an analog of the Gauss-Newton method.
It can be considered as a regularized Gauss-Newton method for solving ill-condition linear algebraic systems.
The advantage of using \eqref{eq30} instead of using (4.1.3) in \cite{R526} is that
one only has to compute the lower upper (LU) decomposition of $A^*A+aI$ once while
the algorithm in \cite{R526} requires computing LU at every step. Note that
computing the LU is the main cost for solving a linear system.
Numerical experiments show that the new method competes favorably with the VR in our experiments.

\end{document}